\title{Combinatorial categories and permutation groups}
\author{Gareth A. Jones\\
School of Mathematics\\
University of Southampton\\
Southampton SO17  1BJ, U.K.\\
{\tt G.A.Jones@maths.soton.ac.uk}
}

\documentclass[10pt]{article}
\usepackage{color}
\usepackage[latin1]{inputenc}
\usepackage{a4wide}
\usepackage{latexsym}
\usepackage{amsmath}
\usepackage{amssymb}
\usepackage{amsfonts}
\newtheorem{thm}{Theorem}[section]
\newtheorem{lemma}[thm]{Lemma}

\date{}
\begin{document}

\maketitle

\begin{abstract}
The regular objects in various categories, such as maps, hypermaps or covering spaces, can be identified with the normal subgroups $N$ of a given group $\Gamma$, with quotient group isomorphic to $\Gamma/N$. It is shown how to enumerate such objects with a given finite automorphism group $G$, how to represent them all as quotients of a single regular object ${\mathcal U}(G)$, and how the outer automorphism group of $\Gamma$ acts on them. Examples constructed include kaleidoscopic maps with trinity symmetry.
\end{abstract}



\section{Introduction}

In certain categories $\mathfrak C$, the objects $\mathcal O$ can be identified with the permutation representations of a particular group $\Gamma=\Gamma_{\mathfrak C}$ on sets $\Phi=\Phi_{\mathcal O}$, and the morphisms ${\mathcal O}\to{\mathcal O}'$ correspond to the functions $\Phi_{\mathcal O}\to\Phi_{{\mathcal O}'}$ commuting with the actions of $\Gamma$. In the case of maps on surfaces one takes $\Gamma$ to be the free product $V_4*C_2$ acting on flags, or $C_{\infty}*C_2$ acting on directed edges of oriented maps. The corresponding groups for hypermaps are $C_2*C_2*C_2$ and the free group $F_2=C_{\infty}*C_{\infty}$ of rank $2$. For abstract polytopes of a given type one can use the corresponding string  Coxeter group, again acting on flags, though here one has to restrict attention to quotient groups satisfying the intersection property. In the case of coverings of a path-connected space $X$ one uses the fundamental group $\pi_1X$, acting on sheets or more precisely on the fibre over a base-point.

In such a case we will call $\mathfrak C$ a permutational category, with parent group $\Gamma$. Each object $\mathcal O$ in such a category $\mathfrak C$ is a disjoint union of connected subobjects, corresponding to the orbits of $\Gamma$ on $\Phi$; one usually restricts attention to the connected objects, as we shall here, so that $\Phi$ can be identified with the set of cosets in $\Gamma$ of a point-stabiliser $M=\Gamma_{\phi}$, where $\phi\in\Phi$. The permutation group induced by $G$ on $\Phi$ is the monodromy group $G={\rm Mon}\,{\mathcal O}={\rm Mon}_{\mathfrak C}{\mathcal O}$ of $\mathcal O$, a subgroup of the symmetric group ${\rm Sym}\,\Phi$ on $\Phi$. The automorphism group $A={\rm Aut}\,{\mathcal O}={\rm Aut}_{\mathfrak C}{\mathcal O}$ of $\mathcal O$, regarded as an object in $\mathfrak C$, is the centraliser of $G$ in ${\rm Sym}\,\Phi$; since $G$ is transitive on $\Phi$, $A$ acts semiregularly on $\Phi$, and
\[A\cong N_{\Gamma}(M)/M\cong N_G(G_{\phi})/G_{\phi}.\] 

The most symmetric objects in $\mathfrak C$ are the regular objects, those for which $A$ acts transitively (and hence regularly) on $\Phi$. This is equivalent to $M$ being a normal subgroup of $\Gamma$, in which case
\[A\cong\Gamma/M\cong G.\]
Indeed, in this case $A$ and $G$ can be identified with the left and right regular representation of the same group.  In principle, understanding regular objects is sufficient for an understanding of all objects in $\mathfrak C$, since each object ${\mathcal O}\in{\mathfrak C}$ is the quotient of some regular object $\tilde{\mathcal O}\in{\mathfrak C}$, corresponding to the core $N$ of $M$ in $\Gamma$, by a group $M/N$ of automorphisms of $\tilde{\mathcal O}$; moreover,  $\tilde{\mathcal O}$ is finite if and only if $\mathcal O$ is finite, since $N$ has finite index in $\Gamma$ if and only if $M$ has finite index. We shall therefore concentrate, for the remainder of this paper, on the regular objects in various categories $\mathfrak C$. In particular, we will study the set ${\mathcal R}(G)={\mathcal R}_{\mathfrak C}(G)$ of regular objects ${\mathcal O}\in{\mathfrak C}$ with ${\rm Aut}\,{\mathcal O}$ isomorphic to a given group $G$. If $\Gamma$ is finitely generated and $G$ is finite then $r(G):=|{\mathcal R}(G)|$ is finite. We will consider how to calculate $r(G)$ in this case, how to represent the objects in ${\mathcal R}(G)$ as quotients of a single regular object ${\mathcal U}(G)$ in $\mathfrak C$, and how the outer automorphism group ${\rm Out}\,\Gamma$ of $\Gamma$ acts on ${\mathcal R}(G)$. Examples will be given, in which the objects are maps, hypermaps or surface coverings, some of them relevant to recent work by Archdeacon, Conder and \v Sir\' a\v n on kaleidoscopic maps with trinity symmetry~\cite{ACS}.

\medskip

\noindent{\bf Acknowledgement} The author thanks the organisers of GEMS 2013 for inviting and supporting him to give a talk summarising this paper. It is based upon work supported by the project: Mobility --- enhancing research, science and education at the Matej Bel University, ITMS code: 26110230082, under the Operational Program Education cofinanced by the European Social Fund.

\section{Examples of permutational categories}

Let us call a category $\mathfrak C$ a {\sl permutational category\/} if it is equivalent to the category of permutation representations of some group $\Gamma$, called the {\sl parent group\/} of $\mathfrak C$. This means that there are functors from each category to the other, so that their composition, in either order, is naturally equivalent to the identity. There are some well-known examples in the literature, though the equivalences are rarely expressed in terms of categories. We will summarise them briefly here; for further details, see, for example,~\cite{JT83} for maps, and~\cite{Jam} for hypermaps.

The category $\mathfrak M$ of maps on surfaces, with branched coverings of maps as its morphisms, is a permutational category, with parent group
\begin{equation}\label{GammaM}
\Gamma=\Gamma_{\mathfrak M}
=\langle R_0, R_1, R_2\mid R_i^2=(R_0R_2)^2=1\rangle.
\end{equation}
Here each $R_i$ acts on the set $\Phi$ of vertex-edge-face flags of a map ${\mathcal O}\in{\mathfrak M}$ by changing, in the only way possible, the $i$-dimensional component of each flag while preserving its $j$-dimensional component for each $j\ne i$. (A boundary flag is fixed by $R_i$ if no such change is possible.) This group, which is a free product
\[\langle R_0, R_2\rangle*\langle R_1\rangle\cong V_4*C_2\]
of a Klein four-group and a cyclic group of order $2$, can be regarded as the extended triangle group $\Delta[\infty, 2, \infty]$ of type $(\infty, 2, \infty)$, generated by reflections in the sides of a hyperbolic triangle with angles $0, \pi/2, 0$. This gives a functor from maps to permutation representations of $\Gamma$. Conversely, given a permutation representation of $\Gamma$ on a set $\Phi$, one can take a set of triangles in bijective correspondence with $\Phi$, each with vertices labelled $0, 1, 2$, and use the cycles of $R_i$ on $\Phi$ to join pairs of triangles across edges $jk\;(j, k\ne i)$; the result is the barycentric subdivision of a map ${\mathcal O}\in{\mathfrak M}$, with the vertices of $\mathcal O$ labelled $0$ and its edges formed by edges of triangles labelled $01$, so that midpoints of edges and faces of $\mathcal O$ are labelled $1$ and $2$. Branched coverings between maps $\mathcal O$ correspond to $\Gamma$-equivariant functions between sets $\Phi$, so we obtain functors ${\mathcal O}\mapsto\Phi$ and $\Phi\mapsto{\mathcal O}$ which give the required equivalence of categories.

Other triangle groups act as parent groups for related categories. For the category ${\mathfrak M}_k$ of maps with all vertex-valencies dividing $k$ we add the relation $(R_1R_2)^k=1$ to the presentation~(\ref{GammaM}), giving the parent group
\begin{equation}\label{GammaMk}
\Gamma_{{\mathfrak M}_k}
=\langle R_0, R_1, R_2\mid R_i^2=(R_0R_2)^2=(R_1R_2)^k=1\rangle
=\Delta[k, 2,\infty].
\end{equation}
Similarly, the isomorphic group $\Delta[\infty,2,k]$ is the parent group for the dual maps, with all face-valencies dividing $k$. For the category $\mathfrak H$ of hypermaps, where hyperedges may be incident with any number of hypervertices and hyperfaces, we delete the relation $(R_0R_2)^2=1$ from~(\ref{GammaM}), giving the parent group
\begin{equation}\label{GammaH}
\Gamma_{\mathfrak H}
=\langle R_0, R_1, R_2\mid R_i^2=1\rangle
=\Delta[\infty,\infty,\infty]
\cong C_2*C_2*C_2
\end{equation}
again permuting flags. Similarly, the extended triangle group
\[\Delta[l,m,n]=\langle R_0, R_1, R_2\mid R_i^2=(R_1R_2)^l=(R_0R_2)^m=(R_0R_1)^n=1\rangle\]
is the parent group for hypermaps of type dividing $(l,m,n)$, that is, of type $(l',m',n')$ where $l', m'$ and $n'$ divide $l, m$ and $n$.

For the corresponding categories ${\mathfrak M}^+$, ${\mathfrak M}_k^+$ and ${\mathfrak H}^+$ of oriented maps and hypermaps we take the orientation-preserving subgroups of index $2$ in these groups, generated by the elements $X=R_1R_0$, $Y=R_0R_2$ and $Z=R_2R_1$ satisfying $XYZ=1$. These are the triangle groups
\begin{equation}\label{GammaM+}
\Gamma_{{\mathfrak M}^+}
=\langle X, Y, Z\mid Y^2=XYZ=1\rangle
=\Delta(\infty,2,\infty)
\cong C_{\infty}*C_2,
\end{equation}
\begin{equation}\label{GammaMk+}
\Gamma_{{\mathfrak M}_k^+}
=\langle X, Y, Z\mid X^k=Y^2=XYZ=1\rangle
=\Delta(k,2,\infty)
\cong C_k*C_2
\end{equation}
and
\begin{equation}\label{GammaH+}
\Gamma_{{\mathfrak H}^+}
=\langle X, Y, Z\mid XYZ=1\rangle
=\Delta(\infty,\infty,\infty)
\cong C_{\infty}*C_{\infty}
\cong F_2.
\end{equation}
Similarly, the triangle group
\[\Delta(l,m,n)=\langle X, Y, Z \mid X^l=Y^m=Z^n=XYZ=1\rangle\]
is the parent group for oriented hypermaps of type dividing $(l,m,n)$.

In the case of oriented hypermaps, the Walsh map~\cite{Wal} represents a hypermap as a bipartite map, with black and white vertices representing hypervertices and hyperedges, and edges representing their incidence; then $X$ and $Y$ permute the set $\Phi$ of edges by following the local orientation around their incident black and white vertices. For oriented maps, $X$ rotates directed edges around their target vertices, while $Y$ reverses them; equivalently, one can convert a map into the Walsh map of a hypermap by adding a white vertex at the centre of each edge, so that new edges correspond to directed edges of the original map.


If $X$ is a path connected, locally path connected, and semilocally simply connected topological space~\cite[Ch.~13]{Mun}), the unbranched coverings $\beta:Y\to X$ of $X$ form a permutational category $\mathfrak C$ with the fundamental group $\Gamma=\pi_1X$ as parent group, using unique path-lifting to permute the fibre $\Phi=\beta^{-1}(x_0)\subset Y$ of $\beta$ over a chosen base-point $x_0\in X$. The regular coverings $\beta$ correspond to the normal subgroups $N$ of $\Gamma$, with covering group ${\rm Aut}\,\beta\cong\Gamma/N$. If $X$ is also a compact Hausdorff space (for instance, a compact manifold or orbifold), then $\Gamma$ is finitely generated~\cite[p.~500]{Mun}.

The categories of maps and hypermaps described above can be regarded as obtained in the above way from suitable orbifolds $X$, such as a triangle with angles $\pi/l,\pi/m, \pi/n$ for hypermaps of type dividing $(l,m,n)$, or a sphere with three cone-points of orders $l,m,n$ in the oriented case. Similarly, Grothendieck's dessins d'enfants~\cite{GG, Gro} are the finite coverings of a sphere minus three points, so their parent group is its fundamental group $\Gamma=F_2$, with generators $X, Y$ and $Z$ inducing the monodromy permutations at the three punctures.

{\sl For the rest of this paper, $\mathfrak C$ will denote a permutationsl category with a finitely generated parent group $\Gamma$.}

\section{Counting regular objects}

For each group $G$, there is a natural bijection between the set ${\mathcal R}(G)={\mathcal R}_{\mathfrak C}(G)$ of (isomorphism classes of) regular objects ${\mathcal O}\in{\mathfrak C}$ with ${\rm Aut}\,{\mathcal O}\cong G$ and the set ${\mathcal N}(G)={\mathcal N}_{\Gamma}(G)$ of normal subgroups $N$ of $\Gamma$ with $\Gamma/N\cong G$. These normal subgroups are the kernels of the epimorphisms $\Gamma\to G$. Two such epimorphisms have the same kernel if and only if they differ by an automorphism of $G$, so there is a bijection between ${\mathcal N}(G)$ and the set of orbits of ${\rm Aut}\,G$, acting by composition on the set ${\rm Epi}(\Gamma,G)$ of epimorphisms $\Gamma\to G$. This action of ${\rm Aut}\,G$ is semiregular, since only the identity automorphism of $G$ fixes an epimorphism.

If $G$ is finite then so is ${\rm Epi}(\Gamma,G)$, since each epimorphism $\Gamma\to G$ is uniquely determined by the images in $G$ of a finite set of generators of $\Gamma$. In this case the sets ${\mathcal R}(G)$ and ${\mathcal N}(G)$ have the same finite cardinality
\begin{equation}\label{r(G)n(G)}
r(G)=r_{\mathfrak C}(G)=|{\mathcal R}(G)|
=n(G)=n_{\Gamma}(G)=|{\mathcal N}(G)|
=\frac{|{\rm Epi}(\Gamma,G)|}{|{\rm Aut}\,G|}.
\end{equation}

In~\cite{Hal}, Hall developed a method for counting epimorphisms onto $G$ by first counting homomorphisms (generally an easier task) to subgroups of $G$, and then using M\"obius inversion in the lattice $\Lambda(G)$ of subgroups of $G$.

Let $\sigma$ and $\phi$ be functions from isomorphism classes of finite groups to $\mathbb C$ such that
\begin{equation}\label{sigma}
\sigma(G)=\sum_{H\le G}\phi(H)
\end{equation}
for all finite groups $G$. Then a simple calculation gives the M\"obius inversion formula for $G$, namely
\begin{equation}\label{phi}
\phi(G)=\sum_{H\le G}\mu_G(H)\sigma(H)
\end{equation}
where $\mu_G$ is the M\"obius function on $\Lambda(G)$, defined recursively by
\begin{equation}\label{mu}
\sum_{K\ge H}\mu_G(K)=\delta_{H,G},
\end{equation}
with $\delta_{H,G}=1$ if $H=G$ and $0$ otherwise. (One can view this as a group-theoretic analogue of the inclusion-exclusion principle, which applies to the lattice of all subsets of $G$; in that situation, by replacing the condition $K\ge H$ in~(\ref{mu}) with $K\supseteq H$ one assigns the value $(-1)^{|G\setminus H|}$ to $\mu_G(H)$ for each subset $H$ of $G$.)

Each homomorphism $\Gamma\to G$ is an epimorphism onto a unique subgroup $H\le G$, so one can take $\sigma(G)$ and $\phi(G)$ to be the numbers of homomorphisms and epimorphisms from $\Gamma$ to $G$ (or possibly those satisfying some extra condition, such as being smooth, i.e.~having a forsion-free kernel). Thus
\begin{equation}\label{Hom}
|{\rm Hom}(\Gamma,G)|=\sum_{H\le G}|{\rm Epi}(\Gamma,H)|,
\end{equation}
so M\"obius inversion gives
\begin{equation}
|{\rm Epi}(\Gamma,G)|=\sum_{H\le G}\mu_G(H)|{\rm Hom}(\Gamma,H)|.
\end{equation}
This proves the first part of the following theorem; the second follows easily.

\begin{thm}\label{enumthm}
If $\mathfrak C$ is a permutational category with a finitely generated parent group $\Gamma$, and $G$ is a finite group, then the number $r(G)$ of isomorphism classes of regular objects ${\mathcal O}\in{\mathfrak C}$ with ${\rm Aut}\,{\mathcal O}\cong G$ is given by
\begin{equation}\label{r(G)}
r(G)=\frac{1}{|{\rm Aut}\,G|}\sum_{H\le G}\mu_G(H)|{\rm Hom}(\Gamma,H)|.
\end{equation}
The number $m(G)$ of isomorphism classes of objects ${\mathcal O}\in{\mathfrak C}$ with ${\rm Mon}\,{\mathcal O}\cong G$ is given by
\begin{equation}\label{m(G)}
m(G)=r(G)c(G),
\end{equation}
where $c(G)$ is the number of conjugacy classes of subgroups of $G$ with trivial core.
\end{thm}

Applying equation~(\ref{r(G)}) to a specific pair $\mathfrak C$ and $G$ requires three ingredients: one must know $|{\rm Aut}\,G|$, $\mu_G(H)$ for each $H\le G$, and $|{\rm Hom}(\Gamma,H)|$ for each $H\le G$ such that $\mu_G(H)\ne 0$. The first is usually the easiest to deal with: for instance  $|{\rm Aut}\,C_n|$ is given by Euler's function $\phi(n)$, while the automorphism groups of the finite simple groups are all known and can be found in sources such as~\cite{ATLAS, Wil}. Finding the other two ingredients is generally more troublesome, and this has been achieved only in special cases.

\section{Evaluating the M\"obius function}

Evaluating the M\"obius function $\mu_G$ requires detailed knowledge of the subgroup lattice of $G$. It has been achieved for several infinite classes of groups $G$, and of course for specific groups which are not too large one can use systems such as GAP or MAGMA. In this context, the database of subgroup lattices described by Connor and Leemans in~\cite{CL1}, and available at~\cite{CL2}, is a valuable resource.

\medskip

\noindent{\bf Example 4.1} A finite cyclic group $G=C_n$ of order $n$ has a unique subgroup $H\cong C_m$ for each $m$ dividing $n$, and no other subgroups. Hall~\cite{Hal} showed that $\mu_G(H)=\mu(n/m)$, where $\mu$ is the M\"obius function of elementary number theory, given by $\mu(n)=(-1)^k$ if $n$ is a product of $k$ distinct primes, and $\mu(n)=0$ otherwise. Indeed, here $\mu$ can be regarded as the M\"obius function on the lattice of subgroups of finite index in the infinite cyclic group $\mathbb Z$.

\medskip

\noindent{\bf Example 4.2} Similarly, it is an easy exercise to compute the M\"obius function for a finite dihedral group; see~\cite{Jon95}.

\medskip

\noindent{\bf Example 4.3} An elementary abelian group $G$ of order $p^d$ can be regarded as a vector space of dimension $d$ over the field ${\mathbb F}_p$, and its subgroups $H$ as the linear subspaces. The number of these of each codimension $k=0,1,\ldots, d$ is equal to the Gaussian binomial coefficient
\[\left(\,\begin{matrix}d \cr k \cr \end{matrix}\,\right)_p
=\frac{(p^d-1)(p^{d-1}-1)\ldots(p^{d-k+1}-1)}{(p^k-1)(p^{k-1}-1)\ldots(p-1)},\]
and Hall~\cite{Hal} showed that they satisfy
\[\mu_G(H)=(-1)^kp^{k(k-1)/2}.\]

Hall showed that in any finite group $G$, if $\mu_G(H)\ne 0$ then $H$ must be the intersection of a set of maximal subgroups of $G$, so in particular $H$ must contain the Frattini subgroup $\Phi(G)$ of $G$, the intersection of all its maximal subgroups.

\medskip

\noindent{\bf Example 4.4}  If $G$ is a $d$-generator finite $p$-group then $\Phi(G)$ is the subgroup $G'G^p$ generated by the commutators and $p$-th powers in $G$, and $G/\Phi(G)$ is an elementary abelian $p$-group of order $p^d$. The subgroups of $H\le G$ with $\mu_G(H)\ne 0$ all contain $\Phi(G)$, and correspond to the subgroups of $G/\Phi(G)$, with $\mu_G(H)$ given by the preceding example.

\medskip

If $G=G_1\times G_2$ where $G_1$ and $G_2$ are finite groups of coprime orders, each subgroup $H\le G$ has the unique form $H=H_1\times H_2$ where $H_i\le G_i$. In this case Hall showed that $\mu_G(H)=\mu_{G_1}(H_1)\mu_{G_2}(H_2)$.

\medskip

\noindent{\bf Example 4.5} Each nilpotent finite group $G$ is a direct product of its Sylow subgroups, which are $p$-groups for the different primes $p$ dividing $|G|$, so the preceding examples show how to compute $\mu_G$.

\medskip

\noindent{\bf Example 4.6} Dickson described the subgroups of the groups $L_2(q)=PSL_2(q)$ in~\cite[Ch.~XII]{Dic}. Using this, Hall~\cite{Hal} calculated the M\"obius function $\mu_G$ for the simple groups $G=L_2(p)$ for primes $p\ge 5$. Equation~(\ref{phi}) takes the form
\[\phi(G)=\sigma(G)-(p+1)\sigma(G_{\infty})-\frac{p(p-1)}{2}\sigma(D_{\frac{p+1}{2}})\]
\[-\frac{p(p+1)}{2}D_{\frac{p-1}{2}}+p(p+1)\sigma(C_{\frac{p-1}{2}})+|G|S,\]
where $G_{\infty}$ is the subgroup of index $p+1$ fixing $\infty$, and $S$ depends on the congruence classes of $p$ mod~$(5)$ and mod~$(8)$, which determine the existence of proper subgroups $H\cong A_5$ or $S_4$. For example, if $p=5$, or if $p\equiv\pm 2$ mod~$(5)$ and $p\equiv\pm 3$ mod~$(8)$, so that there are no such subgroups, then
\[S=-\frac{1}{12}\sigma(A_4)+\frac{1}{4}\sigma(V_4)+\frac{1}{3}\sigma(C_3)+\frac{1}{2}\sigma(C_2)-\sigma(1);\]
there are similar formulae in the other cases. In~\cite{DowPhD} Downs extended Hall's calculation of $\mu_G$ to $L_2(q)$ and $PGL_2(q)$ for all prime powers $q$; see~\cite{DowJLMS} for a proof for $L_2(2^e)$ and a statement of results for $L_2(q)$ where $q$ is odd, and~\cite{DJ} for some combinatorial applications by Downs and the author.

\medskip

\noindent{\bf Example 4.7} The Suzuki groups $G=Sz(q)$ are a family of non-abelian finite simple groups, with $q=2^e$ for some odd $e>1$; see~\cite{ATLAS, Suz, Wil} for their properties, which are similar to those of the groups $L_2(2^e)$. Downs calculated $\mu_G$ for these groups in~\cite{DowSuz}; see~\cite{DJ13} for a statement of the results and some applications.

\section{Counting homomorphisms} 

In order to apply equation~(\ref{r(G)}) to a group $G$, one needs to evaluate $|{\rm Hom}(\Gamma,H)|$ for those subgroups $H\le G$ with $\mu_G(H)\ne 0$. If $\Gamma$ has a presentation with generators $X_i$ and defining relations $R_j$, this is equivalent to counting the solutions $(x_i)$ in $H$ of the equations $R_j(x_i)=1$.

\medskip

\noindent{\bf Example 5.1} If $\Gamma$ is a free product $C_{m_1}*\cdots*C_{m_k}$ of cyclic groups of orders $m_i\in{\mathbb N}\cup\{\infty\}$, then
\[|{\rm Hom}(\Gamma,H)|=\prod_{i=1}^k\sum_{m|m_i}|H|_{m}\]
where $|H|_m$ denotes the number of elements of $H$ of order $m$, and we regard all orders $m$ as dividing $\infty$, so that $\sum_{m|\infty}|H|_{m}=|H|$. For instance, if $\Gamma$ is a free group $F_k$ of rank $k$ then $|{\rm Hom}(\Gamma,H)|=|H|^k$.  Similarly, the torsion theorem for free products~\cite[Theorem~IV.1.6]{LS} implies that a homomorphism $\Gamma\to H$ is smooth if and only if it embeds each finite factor $C_{m_i}$, so the number of such homomorphisms can be found by multiplying $k$ factors equal to $|H|_{m_i}$ or $|H|$ as $m_i$ is finite or infinite.

\medskip

For certain groups $\Gamma$, the character table of $H$ gives $|{\rm Hom}(\Gamma,H)|$.

\medskip

\noindent{\bf Example 5.2} If $\Gamma$ is a polygonal group
\[\Delta(m_1, \ldots, m_k)=\langle X_1, \ldots, X_k\mid X_1^{m_1}=\ldots=X_k^{m_k}=X_1\ldots X_k=1\rangle\]
of type $(m_1, \ldots, m_k)$ for some integers $m_i$, then $|{\rm Hom}(\Gamma, H)|$ can be found by summing the following formula~(\ref{trianglechi}) of Frobenius~\cite{Fro} over all choices of $k$-tuples of conjugacy classes ${\mathcal C}_i$ of elements of orders dividing $m_i$.

\begin{thm}\label{frobchi}
Let ${\mathcal C}_i$ ($i=1, \ldots, k$) be conjugacy classes in a finite group $H$. Then the number of solutions of the equation $x_1\ldots x_k=1$ in $H$, with $x_i\in{\mathcal C}_i$ for $i=1, \ldots, k$, is 
\begin{equation}\label{trianglechi}
\frac{|{\mathcal C}_1|\ldots|{\mathcal C}_k|}{|H|}\sum_{\chi}\frac{\chi(x_1)\ldots\chi(x_k)}{\chi(1)^{k-2}}\end{equation}
where $x_i\in{\mathcal C}_i$ and $\chi$ ranges over the irreducible complex characters of $H$.
\end{thm}
Similarly, the number of smooth homomorphisms $\Gamma\to H$ can be found by restricting the summation to classes of elements of order equal to $m_i$. The case $k=3$ of this theorem, where $\Gamma$ is a triangle group, has often been used in connection with oriented maps and hypermaps: see~\cite{Jon94} and~\cite{JS}, for instance.

\medskip

\noindent{\bf Example 5.3} If $\Gamma$ is an orientable surface group $\Pi_g$, that is, the fundamental group
\[\Pi_g=\pi_1{\mathcal S}_g=\langle A_i, B_i\;(i=1,\ldots, g)\mid\prod_{i=1}^g[A_i,B_i]=1\rangle\]
of a compact orientable surface ${\mathcal S}_g$ of genus $g\ge 1$, one can use the following theorem of Frobenius~\cite{Fro} and Mednykh~\cite{Med}, which counts solutions of the equation $\prod_{i=1}^g[a_i,b_i]=1$:
\begin{thm}\label{FrobMedchi}
If $H$ is any finite group then
\begin{equation}\label{Pi+chi}
|{\rm Hom}(\Pi_g,H)|=|H|^{2g-1}\sum_{\chi}\chi(1)^{2-2g},
\end{equation}
where $\chi$ ranges over the irreducible complex characters of $H$.
\end{thm}

\noindent{\bf Example 5.4} If $\Gamma$ is a non-orientable surface group
\[\Pi_g^{-}=\langle A_i\;(i=1,\ldots, g)\mid\prod_{i=1}^gA_i^2=1\rangle\]
of genus $g\ge 1$, one can use the following result of Frobenius and Schur~\cite{FS}:
\begin{thm}\label{FrobSchchi}
If $H$ is a finite group then
\begin{equation}\label{Pi-chi}
|{\rm Hom}(\Pi^-_g,H)|=|H|^{g-1}\sum_{\chi}c_{\chi}^g\chi(1)^{2-g},
\end{equation}
where $\chi$ ranges over the irreducible complex characters of $H$.
\end{thm}
Here $c_{\chi}$ is the Frobenius-Schur indicator $|H|^{-1}\sum_{h\in H}\chi(h^2)$ of $\chi$, equal to $1, -1$ or $0$ as $\chi$ is the character of a real representation, the real character of a non-real representation, or a non-real character. See~\cite{Jon95} for applications of these two theorems, and~\cite[Ch.~7]{Ser} for several generalisations of them.

\section{Enumerations}

Using Theorem~\ref{enumthm} one can now enumerate, for a given finite group $G$, the regular objects in $\mathfrak C$ with automorphism group $G$, and also the objects in $\mathfrak C$ with monodromy group $G$.

\medskip

\noindent{\bf Example 6.1} It follows from a result of Hall~\cite{Hal} that if $G=L_2(p)$ for some prime $p\ge 5$ and ${\mathfrak C}={\mathfrak H}^+$, so that $\Gamma=F_2$, then
\[r(G)=\frac{1}{4}(p+1)(p^2-2p-1)-\epsilon,\]
where $\epsilon=49$, $40$, $11$ or $2$ as $p\equiv\pm 1$ mod~(5) and $\pm 1$ mod~$(8)$, or $\pm 1$ mod~(5) and $\pm 3$ mod~$(8)$, or $\pm 2$ mod~(5) and $\pm 1$ mod~$(8)$, or $\pm 2$ mod~(5) and $\pm 3$ mod~$(8)$. We also take $\epsilon=2$ when $p=5$, so that $r(G)=19$ in this case; the $19$ regular oriented hypermaps associated with the icosahedral group $G=L_2(5)\cong A_5$ have been described by Breda and the author in~\cite{BJ}. Since this group $G$ has eight conjugacy classes of proper subgroups, all with trivial core since $G$ is simple, it follows from equation~(\ref{m(G)}) that there are $19\times 8=152$ oriented hypermaps with monodromy group $G$, namely the quotients ${\mathcal O}/H$ where ${\mathcal O}\in{\mathcal R}(G)$ and $H<G$.

\medskip

\noindent{\bf Example 6.2} In~\cite{DowPhD}, Downs considered the categories $\mathfrak H$, ${\mathfrak H}^+$, $\mathfrak M$, ${\mathfrak M}^+$, ${\mathfrak M}_3$ and ${\mathfrak M}^+_3$, and gave formulae for $r(G)$ where $G=L_2(q)$ or $PGL_2(q)$ for any prime power $q$. The results for $G=L_2(2^e)$ are given in~\cite{DJ}. Typical results for odd $q$ are:
\[r_{\mathfrak M}(L_2(p^e))=\frac{1}{8e}\sum_{f|e}\mu\left(\frac{e}{f}\right)p^f(p^f-a)\]
for all $p>2$ and odd $e>1$, where $a=2$ or $4$ as $p\equiv 1$ or $-1$ mod~$(4)$, and
\[r_{{\mathfrak M}_3}(PGL_2(p^e))=\frac{3}{4e}\sum_f\mu\left(\frac{e}{f}\right)(p^f-1)\]
for $p>3$ and $e>1$, where the sum is over all factors $f$ of $e$ with $e/f$ odd.

\medskip

\noindent{\bf Example 6.3} Using Downs's calculation of the M\"obius function for $G=Sz(2^e)$ in~\cite{DowSuz}, he and the author have enumerated various combinatorial objects with automorphism group $G$ in~\cite{DJ13}. Typical results are that
\[r_{{\mathfrak H}^+}(G)=\frac{1}{e}\sum_{f|e}\mu\left(\frac{e}{f}\right)2^f(2^{4f}-2^{3f}-9)\]
and
\[r_{\mathfrak M}(G)=\frac{1}{e}\sum_{f|e}\mu\left(\frac{e}{f}\right)(2^f-1)(2^f-2).\]
The second formula, which also gives the number of reflexible maps in ${\mathcal R}_{{\mathfrak M}^+}(G)$, has been obtained by more direct means by Hubard and Leemans in~\cite{HL}.

\medskip

\noindent{\bf Example 6.4} If $G$ is infinite then ${\mathcal R}(G)$ could be finite or infinite. For instance, if ${\mathfrak C}={\mathfrak H}^+$, so that $\Gamma=F_2$, then $r({\mathbb Z}^2)=1$ whereas $r({\mathbb Z})=\aleph_0$.

\section{Universal covers}

For any group $G$, and any $\mathfrak C$, let
\begin{equation}\label{K(G)}
K(G)=K_{\mathfrak C}(G)=\negthinspace\negthinspace\bigcap_{N\in{\mathcal N}(G)}
\negthinspace\negthinspace\negthinspace N.
\end{equation}
This is a normal subgroup of $\Gamma$, so it corresponds to a regular object
\begin{equation}\label{U(G)}
{\mathcal U}(G)={\mathcal U}_{\mathfrak C}(G)=\negthinspace\negthinspace\bigvee_{{\mathcal O}\in{\mathcal R}(G)}
\negthinspace\negthinspace\negthinspace {\mathcal O}
\end{equation}
which we will call the {\it universal cover} for $G$, the smallest object in $\mathfrak C$ covering each ${\mathcal O}\in{\mathcal R}(G)$. This has automorphism group
\begin{equation}\label{overlineG}
\overline G:={\rm Aut}\,{\mathcal U}(G)\cong \Gamma/K(G).
\end{equation}

If $\Gamma$ has generators $X_i\;(i\in I)$ then one can realise $\overline G$ as the subgroup of the cartesian power $G^{{\mathcal R}(G)}$ of $G$ generated by the elements $(x_{i1}, x_{i2},\ldots)$ for $i\in I$, where $x_{ik}$ is the image of $X_i$ in $G={\rm Aut}\,{\mathcal O}_k$ for some numbering ${\mathcal O}_1, {\mathcal O}_2,\ldots$ of the objects ${\mathcal O}_k\in{\mathcal R}(G)$. In particular, $\overline G$ has the same number of generators as $\Gamma$, and it  satisfies all the identical relations satisfied by $G$: for instance, if $G$ is nilpotent of class $c$, is solvable of derived length $d$, or has exponent $e$, then the same applies to $\overline G$. Finally, if $G$ is finite, as we will assume from now on, then so are ${\mathcal U}(G)$ and $\overline G$, with $|\overline G|$ dividing $|G|^r$ where $r=r(G)$.


\medskip

\noindent{\bf Example 7.1} Let ${\mathfrak C}={\mathfrak H}^+$, so that $\Gamma=F_2$. If $G=C_n$ then
$K(G)=\Gamma'\Gamma^n$, so
\[\overline G=\Gamma/\Gamma'\Gamma^n\cong C_n\times C_n.\]
Represented as a bipartite map, the hypermap ${\mathcal U}(G)$ is a regular embedding of the complete bipartite graph $K_{n,n}$ in a surface of genus $(n-1)(n-2)/2$. In fact, we obtain the same universal cover ${\mathcal U}(G)$ and group $\overline G$ whenever $G$ is a $2$-generator abelian group of exponent $n$.

\medskip

This example shows that $\overline G$ can be a rather small subgroup of $G^r$, since $\overline G\cong G^2$ whereas $r>n$. However, if $G$ is a non-abelian finite simple group, then the following result shows that $\overline G=G^r$ for any category $\mathfrak C$; see~\cite{Jon13} for a proof.

\begin{lemma}
Let $N_1,\ldots, N_r$ be distinct normal subgroups of a group $\Gamma$, with each $G_i:=\Gamma/N_i$ non-abelian and simple. If $K=N_1\cap\cdots\cap N_r$ then
\[\Gamma/K\cong G_1\times\cdots\times G_r.\]
\end{lemma}

Taking $\{N_1,\ldots, N_r\}={\mathcal N}_{\Gamma}(G)$, so $G_i\cong G$ for $i=1,\ldots, r$, gives the result.

\medskip

\noindent{\bf Example 7.2}  Let ${\mathfrak C}={\mathfrak H}^+$ again, and let $G=L_2(5)\cong A_5$. By Example~6.1 we have $r(G)=19$, so $\overline G=G^{19}$, of order
\[60^{19}=609359740010496\times 10^{17}\approx 6.1\times 10^{31}.\]

Guralnick and Kantor~\cite{GK} have shown that if $G$ is a non-abelian finite simple group then each non-identity element of $G$ is a member of a generating pair. If such a group $G$ has exponent $e$ then it follows that ${\mathcal U}_{{\mathfrak H}^+}(G)$ has type $(e,e,e)$, so by the Riemann-Hurwitz formula it has genus
\[g=1+\frac{e-3}{2e}|G|^r.\]
In Example~7.2, for instance, $G$ has exponent $30$, so ${\mathcal U}_{{\mathfrak H}^+}(G)$ has genus
\[1+\frac{9}{20}\times 60^{19}=274218830047232000000000000000001\approx2.742\times 10^{31}.\]

\medskip

For any finite group $G$ we have $|{\rm Epi}(F_2,G)|\le|G|^2$, so
\[r_{{\mathfrak H}^+}(G)\le\frac{|G|^2}{|{\rm Aut}\,G|}=\frac{|G|.|Z(G)|}{|{\rm Out}\,G|}\]
where ${\rm Out}\,G$ is the outer automorphism group ${\rm Aut}\,G/{\rm Inn}\,G$ of $G$. In particular, if $G$ has trivial centre then
\begin{equation}\label{upperbound}
r_{{\mathfrak H}^+}(G)\le\frac{|G|}{|{\rm Out}\,G|}.
\end{equation}
If $G$ is a non-abelian finite simple group, then a randomly-chosen pair of elements generate $G$ with probability approaching $1$ as $|G|\to\infty$: this was proved by Dixon~\cite{Dix} for the alternating groups, Kantor and Lubotzky~\cite{KL} for the classical groups of Lie type, and Liebeck and Shalev~\cite{LS} for the exceptional groups of Lie type. Moreover, convergence is quite rapid. It follows that for such groups the upper bound in~(\ref{upperbound}) is asymptotically sharp, that is,
\[r_{{\mathfrak H}^+}(G)\sim\frac{|G|}{|{\rm Out}\,G|}\quad{\rm as}\quad |G|\to\infty.\]
The information in~\cite{ATLAS, Wil} shows that for each of the infinite families of non-abelian finite simple groups, $|{\rm Out}\,G|$ grows much more slowly than $|G|$, so that $r_{{\mathfrak H}^+}(G)$ grows almost as quickly as $|G|$. For instance, $r_{{\mathfrak H}^+}(A_n)\sim n!/4$ as $n\to\infty$. (See~\cite{Dix05, MT} for more precise results concerning generating pairs for $A_n$.)

\medskip                                                                                                                    

\noindent{\bf Example~7.3} If $G$ is the Monster, the largest sporadic simple group, then
\[|G|=2^{46}.3^{20}.5^9.7^6.11^2.13^3.17.19.23.29.31.41.47.59.71\]
\[=808017424794512875886459904961710757005754368000000000\]
\[\approx 8.080\times 10^{53}.\]
Since $|{\rm Out}\,G|=1$ we have $r:=r_{{\mathfrak H}^+}(G)\approx|G|$, so
\[|\overline G|=|G|^r\approx|G|^{|G|}\approx(8.080\times 10^{53})^{8.080\times 10^{53}}\approx 10^{10^{55.639}}.\]
Since $G$ has exponent
\[e=2^5.3^3.5^2.7.11.13.17.19.23.29.31.41.47.59.71\]
\[=1165654792878376600800\approx 1.166\times 10^{21}, \]
the universal cover ${\mathcal U}_{{\mathfrak H}^+}(G)$ has type $(e,e,e)$ and genus approximately $|\overline G|/2$.

\medskip

Similar considerations apply to other categorises $\mathfrak C$, though the universal covers ${\mathcal U}(G)$ and their automorphism groups $\overline G$ are usually rather smaller.

\medskip

\noindent{\bf Example 7.4}  If ${\mathfrak C}={\mathfrak M}^+$ and $G=A_5$ then $r(G)=3$: the orientably regular maps in ${\mathcal R}(G)$ are the icosahedron, the dodecahedron and the great dodecahedron, of types $\{3,5\}$, $\{5,3\}$ and $\{5,5\}$, and of genera $0$, $0$ and $4$. It follows that $\overline G=G^3$, of order $216000$, and that ${\mathcal U}_{{\mathfrak M}^+}(G)$ is a map of type $\{15,15\}$ and genus
\[g=1+\frac{11}{60}\times 60^3=39601.\]
Similarly $r(G)=3$ if ${\mathfrak C}={\mathfrak M}$: the three regular maps in ${\mathcal R}_{\mathfrak M}(G)$ are the non-orientable antipodal quotients of  those in ${\mathcal R}_{{\mathfrak M}^+}(G)$, and the same applies to the universal covers ${\mathcal U}(G)$ in these two categories.

\medskip

\noindent{\bf Example 7.5} It follows from Theorems~\ref{enumthm} and~\ref{FrobMedchi} that there are $2016$ regular coverings of an orientable surface of genus $2$ with covering group $G=A_5$~\cite{Jon95}. They have genus $61$, while ${\mathcal U}(G)$ has genus $1+60^{2016}$ and covering group $G^{2016}$.

\section{Operations on categories}

The automorphisms of the parent group $\Gamma$ of $\mathfrak C$ permute the subgroups of $\Gamma$. Since inner automorphisms leave invariant each conjugacy class of subgroups, there is an induced action of the outer automorphism group
\[\Omega=\Omega_{\mathfrak C}:={\rm Out}\,\Gamma={\rm Aut}\,\Gamma/{\rm Inn}\,\Gamma\]
of $\Gamma$ on isomorphism classes of objects in $\mathfrak C$. Since $\Omega$ preserves normality and quotient groups, it leaves ${\mathcal N}(G)$ and hence ${\mathcal R}(G)$ invariant for each group $G$. Here we will consider, for various categories $\mathfrak C$, the isomorphic actions of $\Omega_{\mathfrak C}$ on these pairs of sets. In some cases, $\Gamma$ decomposes as a free product, possibly with amalgamation, in which case the structure theorems for such groups~\cite[\S 7.2]{LS} often allow $\Omega_{\mathfrak C}$ to be determined explicitly. The case ${\mathfrak C}={\mathfrak M}$, with $\Gamma=V_4*C_2$, was dealt with by Thornton and the author in~\cite{JT83}; other cases considered here are similar, so proofs are omitted.

\subsection{Operations on oriented hypermaps}

In the case ${\mathfrak C}={\mathfrak H}^+$, with $\Gamma=F_2$, James~\cite{Jam} interpreted $\Omega$ as the group of all operations on oriented hypermaps. For any integer $n\ge 1$, the automorphism group of $F_n$ is generated by the elementary Nielsen transformations: permuting the free generators, inverting one of them, and multiplying one of them by another~\cite[Theorem~3.2]{MKS}. When $n=2$ one can identify $\Omega={\rm Out}\,\Gamma$ with $GL_2({\mathbb Z})$ through its faithful induced action on the abelianisation $\Gamma^{\rm ab}=\Gamma/\Gamma'\cong{\mathbb Z}^2$ of $\Gamma$~\cite[Ch.~I, Prop.~4.5]{LSch}.

This group $\Omega$ can be decomposed as a free product with amalgamation as follows (see~\cite[\S 7.2]{CM} for presentations of $\Omega$). If we take the images of $X$ and $Y$ as a basis for $\Gamma^{\rm ab}$, then there is a subgroup $\Sigma\cong S_3\cong D_3$ of $\Omega$, generated by the matrices
\[E=\left(\,\begin{matrix}0&1\cr 1&0\cr \end{matrix}\,\right)\quad {\rm and}\quad
\left(\,\begin{matrix}0&1\cr -1&-1\cr \end{matrix}\,\right)\]
of order $2$ and $3$; this group, which simply permutes the three vertex colours of an oriented hypermap, regarded as a tripartite map by stellating its Walsh map, was introduced by Mach\`\i\/ in~\cite{Mac}. The central involution $-I$ of $\Omega$ reverses the orientation of each hypermap and, together with $\Sigma$, generates a subgroup
\[\Omega_1=\Sigma\times\langle -I\rangle\cong S_3\times C_2\cong D_6\]
of $\Omega$ which preserves the genus of each hypermap and permutes the periods in its type. If a hypermap is represented as a bipartite map, then the matrices
\[\left(\,\begin{matrix}-1&0\cr 0&1\cr \end{matrix}\,\right)\quad{\rm and}\quad
\left(\,\begin{matrix}1&0\cr 0&-1\cr \end{matrix}\,\right)\]
reverse the cyclic order of edges around each black or white vertex, while preserving the order around those of the other colour; they are sometimes called Petrie operations, since they preserve the embedded bipartite graph but replace faces with Petrie polygons (closed zig-zag paths), so the genus may be changed. These two matrices, together with $E$, generate a subgroup $\Omega_2\cong D_4$ such that
\[\Omega_0:=\Omega_1\cap\Omega_2=\langle E, -I\rangle\cong V_4\cong D_2\]
and
\[\Omega=\Omega_1*_{\Omega_0}\Omega_2\cong D_6*_{D_2}D_4.\]
The torsion theorem for free products with amalgamation~\cite[Theorem~IV.2.7]{LS} shows that the operations of finite order are the conjugates of the elements of $\Omega_1\cup\Omega_2$, described by Pinto and the author in~\cite{JPi}.

For any $2$-generator group $G$, the orbits of $\Omega$ on ${\mathcal R}(G)$ correspond to the $T_2$-systems in $G$, that is, the orbits of ${\rm Aut}\,F_2\times{\rm Aut}\,G$ acting by composition on ${\rm Epi}(F_2,G)$ and hence on generating pairs for $G$. It is known~\cite{DunTS, NN} that this action is transitive if $G$ is abelian, whereas Garion and Shalev~\cite{GS} have shown that if $G$ is a non-abelian finite simple group then the number of orbits tends to $\infty$ as $|G|\to\infty$.

\medskip

\noindent{\bf Example 8.1} It follows from work of Neumann and Neumann~\cite{NN} that the $19$ hypermaps in ${\mathcal R}(A_5)$ form two orbits of lengths $9$ and $10$ under $\Omega$, which acts as $S_9\times S_{10}$ on them.
Those hypermaps whose type is a permutation of $(2,5,5)$, $(3,3,5)$ or $(3,5,5)^{-}$ form the first orbit, while those of type a permutation of $(2,3,5)$, $(3,5,5)^+$ or $(5,5,5)$ form the other; here the superscript $+$ or $-$ indicates that the generators of order $5$ are or are not conjugate in $A_5$.

\medskip

This example illustrates a useful result of Nielsen~\cite{Nie}, that when $\Gamma=F_2$ the order of the commutator $[x,y]$ is an invariant of the action of $\Omega$ on ${\mathcal R}(G)$ for any group $G$: ihere the order is $3$ or $5$ for the hypermaps in the two orbits. 

\subsection{Operations on all hypermaps}

When ${\mathfrak C}={\mathfrak H}$ we have $\Gamma=C_2*C_2*C_2$, containing $F_2$ as a characteristic subgroup of index $2$. As shown by James~\cite{Jam} there is again an action of $GL_2({\mathbb Z})$ on hypermaps, as described above, but now extended to all hypermaps. In this case $-I$, induced by conjugation by $R_1$, is in the kernel of the action (since any orientation is now ignored), and there is a faithful action on $\mathfrak H$ of the group
\[{\rm Out}\,\Gamma\cong GL_2({\mathbb Z})/\langle -I\rangle\cong PGL_2({\mathbb Z})\cong S_3*_{C_2}V_4.\]

\subsection{Operations on oriented maps}

When ${\mathfrak C}={\mathfrak M}^+$ we have $\Gamma=C_{\infty}*C_2$, with $\Omega={\rm Out}\,\Gamma\cong V_4$. This group $\Omega$ is generated by vertex-face duality, induced by the automorphism of $\Gamma$ transposing $X$ and $Z$, and orientation-reversal, induced by inverting $X$ and fixing $Y$. These two involutions commute, modulo conjugation by $Y$.

If we restrict to the category ${\mathfrak M}_k^+$ of oriented maps of valency dividing $k$, then $\Gamma=C_k*C_2$, with $\Omega$ isomorphic to the multiplicative group ${\mathbb Z}_k^*$ of units mod~$(k)$ provided $k>2$. The elements of $\Omega$ are the operations $H_j$ defined by Wilson in~\cite{WilOp}, raising the rotation of edges around each vertex to its $j$th power, and induced by automorphisms fixing $Y$ and sending $X$ to $X^j$ for $j\in{\mathbb Z}_k^*$. These operations $H_j$, studied by Nedela and \v Skoviera in~\cite{NS}, preserve the embedded graph, but can change the surface. When $k=5$, for instance, $H_2$ transposes the icosahedron and the great dodecahedron.

\subsection{Operations on all maps}

When ${\mathfrak C}={\mathfrak M}$ we have $\Gamma=V_4*C_2$, with $\Omega={\rm Out}\,\Gamma\cong S_3$
induced by the automorphism group of the free factor $\langle R_0, R_2\rangle\cong V_4$ permuting its three involutions $R_0, R_2$ and $R_0R_2$. As shown by Thornton and the author~\cite{JT83}, this group $\Omega$ is simply an algebraic reinterpretation of  the group of operations on regular maps introduced by Wilson in~\cite{WilOp} (see also~\cite{Lin}). It is generated by the classical duality of maps, which transposes vertices and faces by transposing $R_0$ and $R_2$, and the Petrie duality, which transposes faces and Petrie polygons by transposing $R_2$ and $R_0R_2$; these two operations have a product of order $3$ which acts as a triality operation, cyclically permuting the sets of vertices, faces and Petrie polygons of each map. As noted by Wilson, maps admitting trialities but not dualities seem to be rather rare: Poulton and the author have given some infinite families of examples in~\cite{JPo}.

If we restrict to the category ${\mathfrak M}_k$ of maps of valency dividing $k$, then
\[\Gamma=\Delta[k,2,\infty]=\langle R_0, R_1\rangle*_{\langle R_0\rangle}\langle R_0, R_2\rangle
\cong D_k*_{C_2}D_2,\]
where the amalgamated subgroup $C_2$ is generated by a reflection $R_0$ in each factor. If $k>2$ the automorphisms of $D_k$ fixing $R_0$ form a group isomorphic to ${\mathbb Z}_k^*$, sending $R_0R_1$ to $(R_0R_1)^j$ for any $j\in{\mathbb Z}_k^*$, while those of $D_2$ fixing $R_0$ simply permute $R_2$ and $R_0R_2$. These extend to automorphisms of $\Gamma$ which generate a subgroup  ${\mathbb Z}_k^*\times C_2$ of ${\rm Aut}\,\Gamma$: the first factor induces Wilson's operations $H_j$, and the second factor induces Petrie duality. The structure theorems for free products with amalgamation~\cite[\S 7.2]{LS} show that this subgroup maps onto ${\rm Out}\,\Gamma$. Since $H_{-1}$ is induced by conjugation by $R_0$ we find that
\[\Omega\cong({\mathbb Z}_k^*/\{\pm 1\})\times C_2.\]
When $k=3$, with $\Omega\cong C_2$, we obtain the outer automorphism of the extended modular group $\Gamma=PGL_2({\mathbb Z})$ studied by Thornton and the author in~\cite{JT86}.

\subsection{Operations on surface coverings}

If ${\mathcal S}_g$ is an orientable surface of genus $g\ge 1$, and $\Gamma=\pi_1{\mathcal S}_g$, then by the Baer-Dehn-Nielsen Theorem the group $\Omega={\rm Out}\,\Gamma$ is isomorphic to the extended mapping class group ${\rm Mod}^{\pm}({\mathcal S}_g)$ of ${\mathcal S}_g$, that is, the group of isotopy classes of self-homeomorphisms of ${\mathcal S}_g$ (see~\cite[Ch.~8]{FM}). The mapping class group ${\rm Mod}\,({\mathcal S}_g)$ is the subgroup of index $2$ corresponding to the orientation-preserving self-homeomorphisms; both groups are finitely presented, with ${\rm Mod}\,{\mathcal S}_g$ generated by the Dehn twists~\cite[Ch.~3]{FM}. The induced action of ${\rm Mod}^{\pm}({\mathcal S}_g)$ on coverings of ${\mathcal S}_g$ corresponds to the action of $\Omega$ on permutation representations of $\Gamma$.

\medskip

\noindent{\bf Example 8.2} If $g=1$ then $\Gamma\cong{\mathbb Z}^2$ and $\Omega\cong {\rm Mod}^{\pm}({\mathcal S}_1)\cong GL_2({\mathbb Z})$, with ${\rm Mod}\,({\mathcal S}_1)$ corresponding to $SL_2({\mathbb Z})$. This is generated by the Dehn twists corresponding to the elementary matrices
\[\left(\,\begin{matrix}1&1\cr 0&1\cr \end{matrix}\,\right)\quad {\rm and}\quad
\left(\,\begin{matrix}1&0\cr 1&1\cr \end{matrix}\,\right).\]

\section{Invariance under operations}

Although it is natural to regard the regular objects in $\mathfrak C$ as its most symmetric objects, some of these may have additional `external' symmetries, in the sense that they are invariant (up to isomorphism) under some or all of the operations in $\Omega$. Self-dual maps, such as the tetrahedron, are obvious examples. For any $\mathfrak C$ and $G$ the group $K(G)$ defined in~(\ref{K(G)}) is a characteristic subgroup of $\Gamma$, so the corresponding regular object ${\mathcal U}(G)$ is invariant under $\Omega$. This shows that each object ${\mathcal O}\in{\mathfrak C}$, regular or not, is covered by an $\Omega$-invariant regular object ${\mathcal U}(G)\in{\mathfrak C}$, which is finite if and only if $\mathcal O$ is, and which has automorphism group $\overline G$ where $G={\rm Mon}\,{\mathcal O}$. The smallest $\Omega$-invariant regular object covering $\mathcal O$ can be obtained by restricting the normal subgroups $N$ in~(\ref{K(G)}) to those in the appropriate orbit of $\Omega$ on ${\mathcal N}_{\Gamma}(G)$.

Richter, \v Sir\'a\v n and Wang~\cite{RSW} have shown that for infinitely many $k$ there are regular $k$-valent maps which are invariant under the group of operations
\[\Omega_1:=\Omega_{\mathfrak M}\cong S_3\]
(see also~\cite[Theorem~3]{JT83}), while Archdeacon, Conder and \v Sir\' a\v n~\cite{ACS} have recently constructed infinite families of $k$-valent orientably regular maps invariant under both $\Omega_1$ and the group
\[\Omega_2:=\Omega_{{\mathfrak M}^+_k}\cong {\mathbb Z}_k^*.\]
They call these `kaleidoscopic maps with trinity symmetry'. In both cases, examples of such maps can be constructed as maps ${\mathcal U}_{\mathfrak M}(G)$ for finite groups $G$: for instance, the map denoted by $M_n$ in~\cite[Theorem~2.2]{ACS} has this form where $G$ is a dihedral group of order $4n$, with $K(G)=\Gamma''(\Gamma')^n$ in $\Gamma=\Gamma_{\mathfrak M}\cong V_4*C_2$.

The connection is as follows. For orientably regular maps, invariance under the operation $H_{-1}\in\Omega_2$ is equivalent to reflexibility, so one needs to find normal subgroups of $\Gamma$ which are invariant under the actions of $\Omega_1={\rm Out}\,\Gamma$ (i.e.~which are characteristic subgroups of $\Gamma$) and (for kaleidoscopic maps) of $\Omega_2\cong {\mathbb Z}_k^*$, where $k$ is the valency of the corresponding map. For any quotient $G$ of $\Gamma$, these two groups $\Omega_i$ act by permuting the subgroups in ${\mathcal N}_{\Gamma}(G)$, so they leave invariant their intersection $K(G)$; the map ${\mathcal U}(G)$ corresponding to $K(G)$ is therefore kaleidoscopic with trinity symmetry.

\medskip

\noindent{\bf Example 9.1} Let $G=A_5$, so that the three maps ${\mathcal M}_i\;(i=1,2,3)$ in ${\mathcal R}(G)$ are the antipodal quotients of the icosahedron, the dodecahedron and the great dodecahedron (see Example~7.4); these have types $\{3,5\}_5$, $\{5,3\}_5$ and $\{5,5\}_3$ where the subscript denotes Petrie length, as in~\cite[\S 8.6]{CM}. Their join ${\mathcal U}(G)$ is a non-orientable regular map of type $\{15, 15\}_{15}$ and genus $39602$, with automorphism group $\overline G\cong A_5^3$. The groups $\Omega_1$ and $\Omega_2$ permute the three maps ${\mathcal M}_i$ ($\Omega_1$ transitively, while $\Omega_2\cong{\mathbb Z}_{15}^*\cong C_2\times C_4$ has orbits $\{{\mathcal M}_2\}$ and $\{{\mathcal M}_1, {\mathcal M}_3\}$), so ${\mathcal U}(G)$ is kaleidoscopic with trinity symmetry. (This is the non-orientable example constructed by a different method in~\cite[\S 7]{ACS}.)

\medskip

\noindent{\bf Example 9.2} For an orientable example, we can take $G=A_5\times C_2$, so ${\mathcal U}(G)$ is the join of ${\mathcal U}(A_5)$, described in the preceding example, and ${\mathcal U}(C_2)$, a reflexible map of type $\{2, 2\}_2$ on the sphere corresponding to the derived group $K(C_2)=\Gamma'$ of $\Gamma$. This gives an orientable map of type $\{30,30\}_{30}$ and genus $374401$, which is kaleidoscopic with trinity symmetry and has automorphism group $(A_5\times C_2)^3$.

\medskip

More generally, if $G$ is a non-abelian finite simple group which is a quotient of $\Gamma$ (the only ones which are not are $L_3(q)$, $U_3(q)$, $L_4(2^e)$, $U_4(2^e)$, $A_6$, $A_7$, $M_{11}$, $M_{22}$, $M_{23}$ and $McL$, according to~\cite[Theorem~4.16]{Sir}), these constructions yield a pair of non-orientable and orientable kaleidoscopic maps which have trinity symmetry and have automorphism groups $G^r$ and $G^r\times C_2^3$, where $r=r_{\mathfrak M}(G)$.

\medskip

\noindent{\bf Example 9.3}  If $G$ is the Suzuki group $Sz(8)$, of order $2^6.5.7.13=29120$, then $r=14$ by Example~6.3; the resulting maps have types $\{k,k\}_k$ and $\{2k,2k\}_{2k}$ where $k=455$, the least common multiple of the valencies $5, 7$ and $13$ of the vertices in the $14$ maps in ${\mathcal R}(G)$ (see~\cite{DJ13}). The orientable map has genus
\[1+\frac{29120^{14}\times 2^3}{2}\times\frac{453}{910}=1+29120^{13}\times 57984.\]

\medskip

If only trinity symmetry is required, as in~\cite{RSW}, then smaller examples of this type can generally be found, with $r$ dividing $6$, by replacing ${\mathcal U}(G)$ with the join of an orbit of $\Omega_1$ on ${\mathcal R}_{\mathfrak M}(G)$. For instance, if $G=L_2(p)$ for some prime $p\equiv \pm 1$ mod~$(24)$ one can take $r=1$.

\section{Finiteness} 

Throughout this paper, we have generally assumed that the group $G$ is finite. If it is not, then not only can ${\mathcal R}_{\mathfrak C}(G)$ be infinite, it can even split into infinitely many orbits under the action of $\Omega_{\mathfrak C}$.

\medskip

\noindent{\bf Example 10.1} Let ${\mathfrak C}={\mathfrak H}^+$, so that $\Gamma=F_2$, and let $G=\langle \, x, y\mid x^3=y^2 \, \rangle$, the group $\pi_1(S^3\setminus K)$ of the trefoil knot $K$. This group, isomorphic to the three-string braid group $B_3=\langle a, b\mid aba=bab\rangle$ with $x=ab$ and $y=ab^2$, has centre $Z(G)=\langle x^3\rangle\cong C_{\infty}$, with $G/Z(G)\cong C_3*C_2\cong PSL_2({\mathbb Z})$. Dunwoody and Pietrowski~\cite{DP} have shown that the pairs $x_i=x^{3i+1},\, y_i=y^{2i+1}\;(i\in{\mathbb Z})$ all generate $G$ and lie in different $T_2$-systems. The corresponding normal subgroups $N\in{\mathcal N}_{\Gamma}(G)$, the kernels of the epimorphisms $\Gamma\to G$ given by $X\mapsto x_i, Y\mapsto y_i$, therefore all lie in different orbits of the group $\Omega={\rm Out}\,\Gamma\cong GL_2({\mathbb Z})$, as do the corresponding hypermaps in ${\mathcal R}_{\mathfrak C}(G)$. 

\medskip


\end{document}